\documentclass[11pt]{article}
\usepackage{amscd}
\usepackage{amsfonts}
\usepackage{amsmath}
\usepackage{amssymb}
\usepackage{amsthm}
\usepackage{bbm}
\usepackage{CJK}
\usepackage{fancyhdr}
\usepackage{graphicx}
\usepackage{hyperref}
\usepackage{indentfirst}
\usepackage{latexsym}
\usepackage{mathrsfs}
\usepackage{xypic}

\usepackage[top=1in,bottom=1in,left=1.25in,right=1.25in]{geometry}
\def\<{\langle}
\def\>{\rangle}

\def\o{\otimes}

\date{}
\begin{document}
\renewcommand{\baselinestretch}{1.2}
\renewcommand{\arraystretch}{1.0}
\title{\bf BiHom-Lie superalgebra  structures}
\author{{\bf Shengxiang Wang$^{1}$, Shuangjian Guo$^{2}$
        \footnote { Corresponding author,  E-mail: shuangjianguo@126.com}  }\\
1.~ School of Mathematics and Finance, Chuzhou University,
 Chuzhou 239000,  China \\
2.~ School of Mathematics and Statistics, Guizhou University of Finance and\\ Economics, Guizhou  550025, China}
 \maketitle
\begin{center}
\begin{minipage}{13.cm}

{\bf \begin{center} ABSTRACT \end{center}}
The aim of this paper is to introduce the notion of BiHom-Lie superalgebras.
This class of algebras is a  generalization of both BiHom-Lie algebras and Hom-Lie superalgebras.
In this article, we first present two ways to construct BiHom-Lie superalgebras from BiHom-associative superalgebras
and Hom-Lie superalgebras by Yau's twist principle.
Also, we  explore some general classes of BiHom-Lie admissible superalgebras and
describe all these classes via $G$-BiHom-associative superalgebras, where $G$ is a subgroup of the symmetric group $S_{3}$.
Finally, we discuss the concept of $\beta^{k}$-derivation of BiHom-Lie superalgebras
and prove that the set of all $\beta^{k}$-derivation has a natural
BiHom-Lie superalgebra  structure. \\


 \end{minipage}
 \end{center}
 \normalsize\vskip1cm

\section*{INTRODUCTION}
\def\theequation{0. \arabic{equation}}
\setcounter{equation} {0}

As generalizations of Lie algebras, Hom-Lie algebras were introduced motivated by
applications to physics and to deformations of Lie algebras, especially Lie algebras
of vector fields. The notion of Hom-Lie algebras was firstly introduced by Hartwig,
Larsson and Silvestrov  to describe the structure of certain q-deformations
of the Witt and the Virasoro algebras, see \cite{Aizawa, Chaichian, Hartwig, Hu}.
More precisely, a  Hom-Lie algebras are different from
Lie algebras as the Jacobi identity is replaced by a twisted form using a morphism.
This twisted Jacobi identity is called Hom-Jacobi identity given by
\begin{eqnarray*}
[\alpha(x),[y,z]]+[\alpha(y),[z,x]]+[\alpha(z),[x,y]]=0.
\end{eqnarray*}

The twisting of parts of the defining identities was transferred to other algebraic structures.
In \cite{Makhlouf2008, Makhlouf2009,Makhlouf2010},  Makhlouf and  Silvestrov introduced the notions of
 Hom-associative algebras, Hom-coassociative coalgebras, Hom-bialgebras  and Hom-Hopf algebras.
The original definition of a Hom-bialgebra involved two linear maps, one twisting the
associativity condition and the other one twisting the coassociativity condition.
Later, two directions of study on Hom-bialgebras were developed, one in which the two maps coincide (these are still
called Hom-bialgebras) and another one, started in \cite{Caenepeel2011}, where the two maps are assumed to be
inverse to each other (these are called monoidal Hom-bialgebras).

The main tool for constructing examples of Hom-type algebras is the so-called twisting
principle introduced by Yau for Hom-associative algebras and extended afterwards to other
types of Hom-algebras, see \cite{Yau1, Yau2}.
Later, Yau  \cite{Yau3}
proposed the definition of quasitriangular Hom-Hopf algebras
and showed that each quasitriangular Hom-Hopf algebra yields a solution of Hom-Yang-Baxter equations.
Meanwhile, Yau \cite{Yau6} defined the classical Hom-Yang-Baxter equation  in the same manner and studied Hom-Lie bialgebras.
In fact, the quasi-element of quasitriangular Hom-Lie bialgebras is a solution of classical Hom-Yang-Baxter equation.

A categorical interpretation of Hom-associative algebras has been given by Caenepeel and Goyvaerts in \cite{Caenepeel2011}.
To any monoidal category $C$, they associate a new monoidal category $\widetilde{\mathcal{H}}(C)$ and call it a Hom-category.
 They proved that a Hom-associative algebra is just an  algebra in $\widetilde{\mathcal{H}}(\mathcal{M}_k)$,
  where $\mathcal{M}_k$  is the category of linear spaces over a base field $k$.
  The similar results holds for Hom-coassociative coalgebras and Hom-bialgebras.
 Later, Chen et al. \cite{CWZ2014} studied the quasitriangular structures of  monoidal Hom-Hopf algebras
  and gave an equivalent description via a braided monoidal category of Hom-modules.
Many more properties and structures of Hom-Hopf algebras have been developed,
see  \cite{chen2014, Gohr2010, LIU2014, wang&GUO2014, WCZ2012} and references cited therein.

In \cite{Graziani}, Graziani et al. studied Hom-bialgebras and Hom-Lie algebras in a so-called  group Hom-category
and called them  BiHom-bialgebras and BiHom-Lie algebras.
They defined BiHom-bialgebras using two commuting multiplicative linear maps $\alpha,\beta$,
which unify Hom-bialgebras and monoidal Hom-bialgebras by setting $\alpha=\beta$ and $\alpha=\beta^{-1}$ respectively.
Also they extended the enveloping algebras and representations of Hom-Lie algebras to BiHom-Lie algebras.

In \cite{Ammar2010}, Ammar and Makhlouf introduced the notion of Hom-Lie superalgebras,
they gave a classification of Hom-Lie admissible superalgebras and proved a  graded version of Hartwig-Larsson-Silvestrov Theorem.
Later, Ammar, Makhlouf and Saadaoui \cite{Ammar2013} studied the representation and the cohomology of Hom-Lie superalgebras,
 and calculated the derivations and the second cohomology group of $q$-deformed Witt superalgebra.
In \cite{Cao}, Cao and Luo studied Hom-Lie superalgebra structures on finite-dimensional simple Lie superalgebras,
while Yuan, Sun and Liu considered  Hom-Lie superalgebra structures on infinite-dimensional simple Lie superalgebras in \cite{Yuan1}.

Motivated by these results, we generalize the notion of Hom-Lie superalgebras and  BiHom-Lie algebras to BiHom-Lie superalgebras
and study the structures of BiHom-Lie superalgebras and BiHom-Lie admissible superalgebras.
This paper is organized as follows.

In Section 1, we recall some basic definitions and facts related with BiHom-associative algebras  and BiHom-Lie superalgebras.

In Section 2, we introduce the notion of BiHom-Lie superalgebras and
 show that any BiHom-associative algebra gives rise to a BiHom-Lie superalgebra (see Theorem 2.6).
Meanwhile, we show a method  to construct BiHom-Lie superalgebras from Hom-Lie superalgebras by Yau's twist principle (see Theorem 2.7).

In Section 3, we introduce BiHom-Lie admissible superalgebras and more general
$G$-BiHom-associative superalgebras, where G is a subgroup of the symmetric group $S_3$.
 We show that BiHom-Lie admissible superalgebras are $G$-BiHom-associative superalgebras (see Propositions 3.7).
As a corollary, we obtain a
classification of BiHom-Lie admissible superalgebras using the symmetric group
$S_3$.

In Section 4, we study the $\beta^{k}$-derivation of a BiHom-Lie  superalgebra and prove that
the set of all $\beta^{k}$-derivation of a BiHom-Lie superalgebra forms a BiHom-Lie  superalgebra  (see Propositions 4.4 ).
As an application, we prove that the inner derivation is a $\beta^{k+1}$-derivation (see Propositions  4.5).

\section{Preliminaries}
\def\theequation{\arabic{section}.\arabic{equation}}
\setcounter{equation} {0}

In this section we recall some basic definitions and results related to our paper.
Throughout the paper, all algebraic systems are supposed to be over a field ${k}$.
 Any unexplained definitions and notations can be found in \cite{Graziani} and \cite{Sweedler}.

\smallskip

{\it Definition 2.1.}(\cite{Graziani}) A {\it BiHom-associative algebra} is a 4-tuple $(A,\mu,\alpha,\beta)$,
where $A$ is a ${k}$-linear space, $\alpha: A\rightarrow A$, $\beta: A\rightarrow A$ and
 $\mu: A\o A\rightarrow A$  are linear maps, with notation $\mu(a\o b)=ab$,
 satisfying the following conditions, for all  $a,a',a''\in A$:
\begin{eqnarray*}
&&\alpha\circ\beta=\beta\circ\alpha,\\
&&\alpha(aa')=\alpha(a)\alpha(a'),\beta(aa')=\beta(a)\beta(a'),~~~~~~~~~~~~~~~~~~~~~~~~~~~~~~~~~~~~~~~~~~~~~~~~~~~~~~~\\
&&\alpha(a)(a'a'')=(aa')\beta(a'').
\end{eqnarray*}
And the maps $\alpha, \beta$ are called the structure maps of $A$.
\smallskip

Clearly, a Hom-associative algebra $(A,\mu,\alpha)$ can be regarded as the BiHom-associative
algebra $(A,\mu,\alpha,\alpha)$.
\medskip

{\it Definition 2.2.}(\cite{Graziani}) A {\it BiHom-Lie algebra}
is a 4-tuple $(L,[\cdot,\cdot],\alpha,\beta)$,
where $L$ is a ${k}$-linear space, $\alpha: L\rightarrow L$, $\beta: L\rightarrow L$ and
 $[\cdot,\cdot]: L\o L\rightarrow L$ are linear maps,
 satisfying the following conditions, for all  $a,a',a''\in A$:
\begin{eqnarray*}
&&\alpha\circ\beta=\beta\circ\alpha,\\
&&\alpha[a,a']=[\alpha(a),\alpha(a')],\beta[a,a']=[\beta(a),\beta(a')],\\
&&[\beta(a),\alpha(a')]=-[\beta(a'),\alpha(a)].\\
&&[\beta^{2}(a),[\beta(a'),\alpha(a'')]]+[\beta^{2}(a'),[\beta(a''),\alpha(a)]]+[\beta^{2}(a''),[\beta(a),\alpha(a')]]=0.~~~~~~~~~~~~~~
\end{eqnarray*}

Obviously, a Hom-Lie algebra $(L,[\cdot,\cdot],\alpha)$ is a particular case of a BiHom-Lie
algebra, namely $(L,[\cdot,\cdot],\alpha,\alpha)$.
 Conversely, a BiHom-Lie algebra $(L,[\cdot,\cdot],\alpha,\alpha)$ with bijective $\alpha$
is the Hom-Lie algebra $(L,[\cdot,\cdot],\alpha)$ .

\section{BiHom-associative superalgebras and BiHom-Lie superalgebras}
\def\theequation{\arabic{section}. \arabic{equation}}
\setcounter{equation} {0}

In this section, we  will present the notions of  BiHom-associative superalgebras and BiHom-Lie superalgebras,
and  construct BiHom-Lie superalgebras from  BiHom-associative superalgebras and Hom-Lie superalgebras,
as a generalization of results in \cite{Ammar2010} and \cite{Graziani}.

Now, let $V$ be a linear superspace over $k$ that is a $Z_{2}$-graded linear space with a direct sum
$V = V_{0}\oplus V_{1}$. The elements of $V_{j}, j={0,1},$ are said to be homogenous and of parity $j$. The parity of
a homogeneous element $x$ is denoted by $|x|$.
\medskip

\noindent{\bf Definition 2.1.}
A BiHom-associative superalgebra is a 4-tuple $(A,\mu,\alpha,\beta)$,
where $A$ is a superspace, $\alpha: A\rightarrow A$ and $\beta: A\rightarrow A$ are even homomorphisms,
 $\mu: A\o A\rightarrow A$  is an even bilinear map, with notation $\mu(a\o b)=ab$ satisfying
\begin{eqnarray}
&&\alpha\circ\beta=\beta\circ\alpha,\\
&&\alpha(ab)=\alpha(a)\alpha(b),\beta(ab)=\beta(a)\beta(b),\\
&&\alpha(a)(bc)=(ab)\beta(c),
\end{eqnarray}
for all homogeneous elements $a,b,c\in A$.
\medskip

Let $(A,\mu_{A},\alpha_{A},\beta_{A})$ and $(B,\mu_{B},\alpha_{B},\beta_{B})$ be two BiHom-associative superalgebras,
an even homomorphism $f:A\rightarrow B$ is said to be a morphism of BiHom-associative superalgebras if
$
\alpha_{B}\circ f=f\circ\alpha_{A}, ~\beta_{B}\circ f=f\circ\beta_{A}$
and $f\circ \mu_{A}=\mu_{B}\circ (f\o f).
$
\medskip

\noindent{\bf Remark 2.2.}
Assume that $\beta=\alpha$ in Definition 2.1,
then the BiHom-associative superalgebra $(A,\mu,\alpha,\beta)$ is the Hom-associative superalgebra in \cite{Ammar2010}.
If the part of parity one in $(A,\mu,\alpha,\beta)$ is trivial,
 then it is just the BiHom-associative algebra in \cite{Graziani}.
\medskip

\noindent{\bf Definition 2.3.}
A BiHom-Lie superalgebra is a 4-tuple $(L,[\cdot,\cdot],\alpha,\beta)$,
where $L$ is a superspace, $\alpha: L\rightarrow L$ and $\beta: L\rightarrow L$ are even homomorphisms,
 $[\cdot,\cdot]: L\o L\rightarrow L$  is an even bilinear map satisfying
\begin{eqnarray}
&&\alpha\circ\beta=\beta\circ\alpha,\\
&&\alpha[x,y]=[\alpha(x),\alpha(y)],\beta[x,y]=[\beta(x),\beta(y)],\\
&&[\beta(x),\alpha(y)]=-(-1)^{|x||y|}[\beta(y),\alpha(x)].\\
&&\circlearrowleft_{x,y,z}(-1)^{|x||z|}[\beta^{2}(x),[\beta(y),\alpha(z)]]=0,
\end{eqnarray}
for all homogeneous elements $x,y,z\in L$.

Let $(L,[\cdot,\cdot],\alpha,\beta)$ and $(L',[\cdot,\cdot]',\alpha',\beta')$ be two BiHom-Lie superalgebras,
an even homomorphism $f:L\rightarrow L'$ is said to be a morphism of BiHom-Lie superalgebras if
$
\alpha'\circ f=f\circ\alpha, \beta'\circ f=f\circ\beta$
and $f\circ [\cdot,\cdot]=[\cdot,\cdot]'\circ (f\o f).
$
\medskip

\noindent{\bf Example 2.4.}
Let $L=L_{0}\oplus L_{1}$ be a 2-dimensional superspace,
$L_{0}$ is generated by $x$ and $L_{1}$ is generated by $y$  such that $[x,y]=0$.
Then for any commutative even homomorphism $\alpha,\beta: L\rightarrow L$,
 $(L,[\cdot,\cdot],\alpha,\beta)$  is a  BiHom-Lie superalgebra.
\medskip

\noindent{\bf Example 2.5.}
Let $L=L_{0}\oplus L_{1}$ be a 3-dimensional superspace,
$L_{0}$ is generated by $e_{1},e_{2}$ and $L_{1}$ is generated by $e_{3}$.
Define a bracket product $[\cdot,\cdot]$ on $L$ by
 $$[e_{1},e_{2}]=e_{1}, [e_{1},e_{3}]=[e_{2},e_{3}]=[e_{3},e_{3}]=0.$$
Let $\lambda,\mu$ be two nonzero scalars in $k$.
Consider the maps $\alpha,\beta: L\rightarrow L$ defined on the basis elements by
 \begin{eqnarray*}
 &&\alpha(e_{1})=\mu(e_{1}), \alpha(e_{2})=e_{2},\alpha(e_{3})=\lambda e_{3},\\
&&\beta(e_{1})=\mu(e_{1}), \beta(e_{2})=e_{2},\beta(e_{3})=-\lambda e_{3}.
  \end{eqnarray*}
It is straightforward to check that  $\alpha,\beta$ defines two BiHom-Lie superalgebra homomorphisms
and $\alpha\circ\beta=\beta\circ\alpha$.
Also one may check that the bracket product $[\cdot,\cdot]$  and the structure maps $\alpha,\beta$
satisfy Eq. (2.6) and Eq. (2.7), then $(L,[\cdot,\cdot],\alpha,\beta)$  is a  BiHom-Lie superalgebra.
\smallskip

\noindent{\bf Theorem 2.6.}
Let $(A,\mu,\alpha,\beta)$ be a BiHom-associative superalgebra with bijective homomorphisms $\alpha$ and $\beta$.
One can define the supercommutator on homogeneous elements by
 \begin{eqnarray*}
[x,y]=xy-(-1)^{|x||y|}\alpha^{-1}(\beta(y))\alpha(\beta^{-1}(x))
  \end{eqnarray*}
and then extending by linearity to all elements. Then $(A, [\cdot,\cdot],\alpha,\beta)$ is a  BiHom-Lie superalgebra.
\smallskip

{\bf Proof}
First we check that the bracket product $[\cdot,\cdot]$ is compatible with the structure maps $\alpha$ and $\beta$.
For any homogeneous elements $x,y\in A$, we have
 \begin{eqnarray*}
[\alpha(x),\alpha(y)]
&=&\alpha(x)\alpha(y)-(-1)^{|\alpha(x)||\alpha(y)|}\alpha^{-1}(\beta(\alpha(y)))\alpha(\beta^{-1}(\alpha(x)))\\
&=&\alpha(x)\alpha(y)-(-1)^{|x||y|}\beta(y)\alpha^{2}(\beta^{-1}(x))\\
&=&\alpha[x,y].
  \end{eqnarray*}
The second equality holds since $\alpha$ is even and $\alpha\circ\beta=\beta\circ\alpha$.
Similarly, one can prove that $\beta[x,y]=[\beta(x),\beta(y)]$.

To verify the skew-supersymmetry, let $x,y\in A$. Then
\begin{eqnarray*}
[\beta(x),\alpha(y)]
&=&\beta(x)\alpha(y)-(-1)^{|\beta(x)||\alpha(y)|}\alpha^{-1}(\beta(\alpha(y)))\alpha(\beta^{-1}(\beta(x)))\\
&=&\beta(x)\alpha(y)-(-1)^{|x||y|}\beta(y)\alpha(x).
  \end{eqnarray*}
Similarly, $[\beta(x),\alpha(y)]=\beta(y)\alpha(x)-(-1)^{|y||x|}\beta(x)\alpha(y)=-(-1)^{|y||x|}[\beta(x),\alpha(y)]$.
So Eq. (2.6) holds.

Now we prove the Eq. (2.7). For any $x,y,z\in A$, we have
\begin{eqnarray*}
&&(-1)^{|x||z|}[\beta^{2}(x),[\beta(y),\alpha(z)]]\\
&&~~~~~~=(-1)^{|x||z|}[\beta^{2}(x),\beta(y)\alpha(z)-(-1)^{|y||z|}\alpha^{-1}(\beta(\alpha(z)))\alpha(\beta^{-1}(\beta(y)))]\\
&&~~~~~~=(-1)^{|x||z|}[\beta^{2}(x),\beta(y)\alpha(z)-(-1)^{|y||z|}\beta(z)\alpha(y)]\\
&&~~~~~~=(-1)^{|x||z|}\beta^{2}(x)(\beta(y)\alpha(z))-(-1)^{|x||y|}(\alpha^{-1}(\beta^{2}(y))\beta(z))\alpha(\beta(x))\\
&&~~~~~~~~~-(-1)^{|x||z|+|y||z|}\beta^{2}(x)(\beta(z)\alpha(y))+(-1)^{|z||y|+|x||y|}(\alpha^{-1}(\beta^{2}(z))\beta(y))\alpha(\beta(x)).
  \end{eqnarray*}
Similarly, we have
\begin{eqnarray*}
&&(-1)^{|y||x|}[\beta^{2}(y),[\beta(z),\alpha(x)]]\\
&&~~~~~~=(-1)^{|y||x|}\beta^{2}(y)(\beta(z)\alpha(x))-(-1)^{|y||z|}(\alpha^{-1}(\beta^{2}(z))\beta(x))\alpha(\beta(y))\\
&&~~~~~~~~~-(-1)^{|y||x|+|z||x|}\beta^{2}(y)(\beta(x)\alpha(z))+(-1)^{|x||z|+|y||z|}(\alpha^{-1}(\beta^{2}(x))\beta(z))\alpha(\beta(y)),\\
&&(-1)^{|z||y|}[\beta^{2}(z),[\beta(x),\alpha(y)]]\\
&&~~~~~~=(-1)^{|z||y|}\beta^{2}(z)(\beta(x)\alpha(y))-(-1)^{|z||x|}(\alpha^{-1}(\beta^{2}(x))\beta(y))\alpha(\beta(z))\\
&&~~~~~~~~~-(-1)^{|z||y|+|x||y|}\beta^{2}(z)(\beta(y)\alpha(x))+(-1)^{|y||x|+|z||x|}(\alpha^{-1}(\beta^{2}(y))\beta(x))\alpha(\beta(z)).
  \end{eqnarray*}
By the associativity Eq. (2.3), it is not hard to check that
\begin{eqnarray*}
&&\circlearrowleft_{x,y,z}(-1)^{|x||z|}[\beta^{2}(x),[\beta(y),\alpha(z)]]=0,
\end{eqnarray*}
as desired. And this finishes the proof.  \hfill $\square$
\medskip

\noindent{\bf Theorem 2.7.}
Let $(L,[\cdot,\cdot])$ be a Lie superalgebra.
Assume that $\alpha,\beta$ are two even commuting algebra homomorphisms of $L$.
Then  $(L,[\cdot,\cdot]_{\alpha,\beta},\alpha,\beta)$, where
$
[x,y]_{\alpha,\beta}=[\alpha(x),\beta(y)],
$
 is a  BiHom-Lie superalgebra.
\smallskip

{\bf Proof}
For any $x,y\in L$, we have
\begin{eqnarray*}
&&[\beta(x),\alpha(y)]_{\alpha,\beta}=[\alpha\beta(x),\beta\alpha(y)]=\alpha\beta([x,y]),\\
&&[\beta(y),\alpha(x)]_{\alpha,\beta}=[\alpha\beta(y),\beta\alpha(x)]=\alpha\beta([y,x])=(-1)^{|x||y|}\alpha\beta([x,y]).
\end{eqnarray*}
So $[\beta(x),\alpha(y)]_{\alpha,\beta}=(-1)^{|x||y|}[\beta(y),\alpha(x)]_{\alpha,\beta}$,
that is, Eq. (2.6) holds.

For Eq. (2.7), we have
\begin{eqnarray*}
&&\circlearrowleft_{x,y,z}(-1)^{|x||z|}[\beta^{2}(x),[\beta(y),\alpha(z)]_{\alpha,\beta}]_{\alpha,\beta}\\
&&~~=\circlearrowleft_{x,y,z}(-1)^{|x||z|}[\beta^{2}(x),[\alpha\beta(y),\alpha\beta(z)]]_{\alpha,\beta}\\
&&~~=\circlearrowleft_{x,y,z}(-1)^{|x||z|}[\alpha\beta^{2}(x),[\alpha\beta^{2}(y),\alpha\beta^{2}(z)]]=0.
\end{eqnarray*}
The last equality holds since $(L,[\cdot,\cdot])$ is a Lie superalgebra.
Thus $(L,[\cdot,\cdot]_{\alpha,\beta},\alpha,\beta)$  is a  BiHom-Lie superalgebra.
\hfill $\square$

\section{BiHom-Lie admissible superalgebras}
\def\theequation{\arabic{section}. \arabic{equation}}
\setcounter{equation} {0}

 In this section, we introduce the notion of BiHom-Lie admissible superalgebras and
provide a classification of BiHom-Lie admissible superalgebras using the symmetric group $S_{3}$.
In this section, we always assume that the structure maps $\alpha$ and $\beta$ are bijective.
\medskip

A BiHom-superalgebra is a 4-tuple $(V,\mu,\alpha,\beta)$,
where $V$ is a superspace, $\alpha: V\rightarrow V$ and $\beta: V\rightarrow V$ are even homomorphism,
 $\mu: V\o V\rightarrow V$  is an even bilinear map satisfying
\begin{eqnarray*}
\alpha\circ\beta=\beta\circ\alpha,
\alpha\circ\mu=\mu\circ(\alpha\o\alpha),\beta\circ\mu=\mu\circ(\beta\o\beta).
\end{eqnarray*}


\noindent{\bf Definition 3.1.}
Let $A=(V,\mu,\alpha,\beta)$ be a BiHom-superalgebra.
Then $A$ is said to be a BiHom-Lie admissible superalgebra over $V$
if the bracket defined by
\begin{eqnarray}
[x,y]=\mu(x\o y)-(-1)^{|x||y|}\mu(\alpha^{-1}(\beta(y))\o\alpha(\beta^{-1}(x)))
\end{eqnarray}
satisfies the BiHom-superJacobi identity (2.7),
for all homogeneous elements $x,y\in V$.
\medskip

\noindent{\bf Remark 3.2.}
By Theorem 2.5, any BiHom-associative superalgebra is a BiHom-Lie admissible superalgebra.
\medskip

Let $(L,[\cdot,\cdot],\alpha,\beta)$ be a BiHom-Lie  superalgebra.
Define a new supercommutator bracket $[\cdot,\cdot]'$  on $L$ by
\begin{eqnarray*}
[x,y]'=[x,y]-(-1)^{|x||y|}[\alpha^{-1}(\beta(y)),\alpha(\beta^{-1}(x))].
\end{eqnarray*}
It is easy to see that the bracket $[\cdot,\cdot]'$ satisfies Eq. (2.6).
Moreover, we have
\begin{eqnarray*}
&&(-1)^{|x||z|}[\beta^{2}(x),[\beta(y),\alpha(z)]']'\\
&&~~~~~~=(-1)^{|x||z|}[\beta^{2}(x),[\beta(y),\alpha(z)]-(-1)^{|y||z|}[\alpha^{-1}(\beta(\alpha(z))),\alpha(\beta^{-1}(\beta(y)))]]'\\
&&~~~~~~=(-1)^{|x||z|}[\beta^{2}(x),[\beta(y),\alpha(z)]-(-1)^{|y||z|}[\beta(z),\alpha(y)]]'\\
&&~~~~~~=(-1)^{|x||z|}[\beta^{2}(x),2[\beta(y),\alpha(z)]]\\
&&~~~~~~=2(-1)^{|x||z|}[\beta^{2}(x),[\beta(y),\alpha(z)]]-2(-1)^{|x||y|}[\alpha^{-1}\beta([\beta(y),\alpha(z)]),\alpha(\beta(x))]\\
&&~~~~~~=4(-1)^{|x||z|}[\beta^{2}(x),[\beta(y),\alpha(z)]].
\end{eqnarray*}
Therefore,
\begin{eqnarray*}
&&\circlearrowleft_{x,y,z}(-1)^{|x||z|}[\beta^{2}(x),[\beta(y),\alpha(z)]']'
=4\circlearrowleft_{x,y,z}(-1)^{|x||z|}[\beta^{2}(x),[\beta(y),\alpha(z)]]
=0.
\end{eqnarray*}

Our discussion above shows:
\medskip

 \noindent{\bf Proposition 3.3.}
Any BiHom-Lie superalgebra is  a  BiHom-Lie  admissible superalgebra.
\smallskip

Let $A=(V,\mu,\alpha,\beta)$ be a BiHom-superalgebra.
The $(\alpha,\beta)$-associator of the  multiplication $\mu$ is a trilinear map $\textbf{as}_{\alpha,\beta}$ on $V$
defined by
\begin{eqnarray*}
\textbf{as}_{\alpha,\beta}(x_{1},x_{2},x_{3})
=\mu(\alpha(x_{1}),\mu(x_{2},x_{3}))-\mu(\mu(x_{1},x_{2}),\beta(x_{3})),
\end{eqnarray*}
where $x_{1},x_{2},x_{3}$ are homogeneous elements in $V$.
\medskip

Now let us introduce the notation:
\begin{eqnarray*}
S(x,y,z):=
\circlearrowleft_{x,y,z}(-1)^{|x||z|}\textbf{as}_{\alpha,\beta}(\alpha^{-1}\beta^{2}(x),\beta(y),\alpha(z)).
\end{eqnarray*}
Then we have the following lemmas:
\medskip

 \noindent{\bf Lemma 3.4.}
Let $A=(V,\mu,\alpha,\beta)$ be a BiHom-superalgebra and $[\cdot,\cdot]$ the associated supercommutator.
Then
\begin{eqnarray*}
&&\circlearrowleft_{x,y,z}(-1)^{|x||z|}[\beta^{2}(x),[\beta(y),\alpha(z)]]
=(-1)^{|x||z|}\textbf{as}_{\alpha,\beta}(\alpha^{-1}\beta^{2}(x),\beta(y),\alpha(z))\\
&&~~+(-1)^{|x||y|}\textbf{as}_{\alpha,\beta}(\alpha^{-1}\beta^{2}(y),\beta(z),\alpha(x))
+(-1)^{|z||y|}\textbf{as}_{\alpha,\beta}(\alpha^{-1}\beta^{2}(z),\beta(x),\alpha(y))\\
 &&~~-(-1)^{|x||z|+|z||y|}\textbf{as}_{\alpha,\beta}(\alpha^{-1}\beta^{2}(x),\beta(z),\alpha(y))
-(-1)^{|x||y|+|z||y|}\textbf{as}_{\alpha,\beta}(\alpha^{-1}\beta^{2}(z),\beta(y),\alpha(x))\\
&&~~-(-1)^{|x||y|+|z||x|}\textbf{as}_{\alpha,\beta}(\alpha^{-1}\beta^{2}(y),\beta(x),\alpha(z)).
\end{eqnarray*}

{\bf Proof} For any homogeneous elements $x,y,z\in V$, we have
\begin{eqnarray*}
&&\circlearrowleft_{x,y,z}(-1)^{|x||z|}[\beta^{2}(x),[\beta(y),\alpha(z)]]\\
&&=(-1)^{|x||z|}\{\mu(\beta^{2}(x)\o\mu(\beta(y)\o\alpha(z)))-\mu(\mu(\alpha^{-1}(\beta^{2}(x))\o\beta(y))\o\alpha(\beta(z)))\}\\
&&~~+(-1)^{|x||y|}\{\mu(\beta^{2}(y)\o\mu(\beta(z)\o\alpha(x)))-\mu(\mu(\alpha^{-1}(\beta^{2}(y))\o\beta(z))\o\alpha(\beta(x)))\}\\
&&~~+(-1)^{|z||y|}\{\mu(\beta^{2}(z)\o\mu(\beta(x)\o\alpha(y)))-\mu(\mu(\alpha^{-1}(\beta^{2}(z))\o\beta(x))\o\alpha(\beta(y)))\}\\
&&~~-(-1)^{|x||z|+|z||y|}\{\mu(\mu(\alpha^{-1}\beta^{2}(x)\o\beta(z))\o\alpha\beta(y))-\mu(\beta^{2}(x)\o\mu(\beta(z)\o\alpha(y))\}\\
&&~~-(-1)^{|x||y|+|z||y|}\{\mu(\mu(\alpha^{-1}\beta^{2}(z)\o\beta(y))\o\alpha\beta(x))-\mu(\beta^{2}(z)\o\mu(\beta(y)\o\alpha(x))\}\\
&&~~-(-1)^{|x||y|+|z||x|}\{\mu(\mu(\alpha^{-1}\beta^{2}(y)\o\beta(x))\o\alpha\beta(z))-\mu(\beta^{2}(y)\o\mu(\beta(x)\o\alpha(z))\}\\
&&=(-1)^{|x||z|}\textbf{as}_{\alpha,\beta}(\alpha^{-1}\beta^{2}(x),\beta(y),\alpha(z))
     +(-1)^{|x||y|}\textbf{as}_{\alpha,\beta}(\alpha^{-1}\beta^{2}(y),\beta(z),\alpha(x))\\
&&~~+(-1)^{|z||y|}\textbf{as}_{\alpha,\beta}(\alpha^{-1}\beta^{2}(z),\beta(x),\alpha(y))
     -(-1)^{|x||z|+|z||y|}\textbf{as}_{\alpha,\beta}(\alpha^{-1}\beta^{2}(x),\beta(z),\alpha(y))\\
&&~~-(-1)^{|x||y|+|z||y|}\textbf{as}_{\alpha,\beta}(\alpha^{-1}\beta^{2}(z),\beta(y),\alpha(x))\\
&&~~-(-1)^{|x||y|+|z||x|}\textbf{as}_{\alpha,\beta}(\alpha^{-1}\beta^{2}(y),\beta(x),\alpha(z)).
\end{eqnarray*}

\medskip

 \noindent{\bf Proposition 3.5.}
Let $A=(V,\mu,\alpha,\beta)$ be a BiHom-superalgebra.
Then $A$ is a BiHom-Lie admissible superalgebra if and only if it satisfies
\begin{eqnarray*}
S(x,y,z)=(-1)^{|x||z|+|z||x|+|x||y|}S(x,z,y),
\end{eqnarray*}
for all homogeneous elements $x,y,z\in V$.
\smallskip

{\bf Proof} For any homogeneous elements $x,y,z\in V$, it is easy to check that
\begin{eqnarray*}
(-1)^{|x||z|+|z||x|+|x||y|}S(x,z,y)
&=&(-1)^{|x||z|+|z||y|}\textbf{as}_{\alpha,\beta}(\alpha^{-1}\beta^{2}(x),\beta(z),\alpha(y))\\
&&+(-1)^{|x||y|+|z||y|}\textbf{as}_{\alpha,\beta}(\alpha^{-1}\beta^{2}(z),\beta(y),\alpha(x))\\
&&+(-1)^{|x||y|+|z||x|}\textbf{as}_{\alpha,\beta}(\alpha^{-1}\beta^{2}(y),\beta(x),\alpha(z)).
\end{eqnarray*}
Therefore, by Lemma 3.4, we have
\begin{eqnarray*}
&&S(x,y,z)-(-1)^{|x||z|+|z||x|+|x||y|}S(x,z,y)\\
&&~~=(-1)^{|x||z|}\textbf{as}_{\alpha,\beta}(\alpha^{-1}\beta^{2}(x),\beta(y),\alpha(z))
     +(-1)^{|x||y|}\textbf{as}_{\alpha,\beta}(\alpha^{-1}\beta^{2}(y),\beta(z),\alpha(x))\\
&&~~+(-1)^{|z||y|}\textbf{as}_{\alpha,\beta}(\alpha^{-1}\beta^{2}(z),\beta(x),\alpha(y))
     -(-1)^{|x||z|+|z||y|}\textbf{as}_{\alpha,\beta}(\alpha^{-1}\beta^{2}(x),\beta(z),\alpha(y))\\
&&~~-(-1)^{|x||y|+|z||y|}\textbf{as}_{\alpha,\beta}(\alpha^{-1}\beta^{2}(z),\beta(y),\alpha(x))\\
&&~~-(-1)^{|x||y|+|z||x|}\textbf{as}_{\alpha,\beta}(\alpha^{-1}\beta^{2}(y),\beta(x),\alpha(z))\\
&&~~=\circlearrowleft_{x,y,z}(-1)^{|x||z|}[\beta^{2}(x),[\beta(y),\alpha(z)]]
\end{eqnarray*}
So $\circlearrowleft_{x,y,z}(-1)^{|x||z|}[\beta^{2}(x),[\beta(y),\alpha(z)]]=0$
if and only if it satisfies
$$
S(x,y,z)-(-1)^{|x||z|+|z||x|+|x||y|}S(x,z,y)=0.
$$
The proof is completed.
\hfill $\square$
\medskip

In the following, we will provide a classification of BiHom-Lie admissible superalgebras using the symmetric group $S_{3}$,
whereas it was classified in \cite{Ammar2010, Makhlouf2008, Yuan} for Hom-Lie admissible algebras,  Hom-Lie admissible superalgebras
and Hom-Lie color admissible algebras, respectively.
\smallskip

Let $S_{3}$ be the symmetric group generated by $\sigma_{1}=(12),\sigma_{2}=(23)$ and
 $A=(V,\mu,\alpha,\beta)$ a BiHom-superalgebra.
Suppose that $S_{3}$ acts on  $V^{\times 3}$ in the usual way, i.e.,
$\sigma(x_{1},x_{2},x_{3})=(x_{\sigma(1)},x_{\sigma(2)},x_{\sigma(3)})$.

For convenience, define the parity of the transposition  $\sigma_{i}$ with $i\in\{1,2\}$ as follows:
\begin{eqnarray*}
|\sigma_{i}(x_{1},x_{2},x_{3})|=|x_{i}||x_{i+1}|.
\end{eqnarray*}
It is natural to assume that the parity of the identity is 0 and for the composition $\sigma_{i}\sigma_{j}$,
it is defined by
\begin{eqnarray*}
|\sigma_{i}\sigma_{j}(x_{1},x_{2},x_{3})|
&=&|\sigma_{j}(x_{1},x_{2},x_{3})|+|\sigma_{i}(\sigma_{j}(x_{1},x_{2},x_{3}))|\\
&=&|\sigma_{j}(x_{1},x_{2},x_{3})|+|\sigma_{i}(x_{\sigma_{j}(1)},x_{\sigma_{j}(2)},x_{\sigma_{j}(3)})|.
\end{eqnarray*}
One can define by induction the parity for any composition as follows:
\begin{eqnarray*}
&&|id(x_{1},x_{2},x_{3})|=0,\\
&&|\sigma_{1}(x_{1},x_{2},x_{3})|=|x_{1}||x_{2}|,\\
&&|\sigma_{2}(x_{1},x_{2},x_{3})|=|x_{2}||x_{3}|,\\
&&|\sigma_{1}\sigma_{2}(x_{1},x_{2},x_{3})|=|x_{2}||x_{3}|+|x_{1}||x_{3}|,\\
&&|\sigma_{2}\sigma_{1}(x_{1},x_{2},x_{3})|=|x_{1}||x_{2}|+|x_{1}||x_{3}|,\\
&&|\sigma_{2}\sigma_{1}\sigma_{2}(x_{1},x_{2},x_{3})|=|x_{2}||x_{3}|+|x_{1}||x_{3}|+|x_{1}||x_{2}|,
\end{eqnarray*}
where $x_{1},x_{2},x_{3}$ are homogeneous element in $V$.
\medskip

 \noindent{\bf Lemma 3.6.}
 A BiHom-superalgebra $A=(V,\mu,\alpha,\beta)$ is a BiHom-Lie admisssible superalgebra
 if and only if the following condition holds
\begin{eqnarray*}
\sum_{\sigma\in S_{3}}(-1)^{\varepsilon(\sigma)}(-1)^{|\sigma(x_{1},x_{2},x_{3})|}
    \textbf{as}_{\alpha,\beta}\circ\sigma(\alpha^{-1}\beta^{2}(x_{1}),\beta(x_{2}),\alpha(x_{3}))=0,
\end{eqnarray*}
for all homogeneous elements $x_{1},x_{2},x_{3}\in V$, where $(-1)^{\varepsilon(\sigma)}$
is the signature of  $\sigma$.
\smallskip

{\bf Proof}
It is sufficient to verify the BiHom-superJacobi identity (2.7).
By Lemma 3.4,
\begin{eqnarray*}
&&\circlearrowleft_{x_{1},x_{2},x_{3}}(-1)^{|x_{1}||x_{3}|}[\beta^{2}(x_{1}),[\beta(x_{2}),\alpha(x_{3})]]\\
&&~~~~=(-1)^{|x_{1}||x_{3}|}\sum_{\sigma\in S_{3}}(-1)^{\varepsilon(\sigma)}(-1)^{|\sigma(x_{1},x_{2},x_{3})|}
    \textbf{as}_{\alpha,\beta}\circ\sigma(\alpha^{-1}\beta^{2}(x_{1}),\beta(x_{2}),\alpha(x_{3})),
\end{eqnarray*}
since $\alpha,\beta$ are even homomorphism.
\hfill $\square$
\medskip

Let $G$ be a subgroup of $S_3$, any BiHom-superalgebra $(V,\mu,\alpha,\beta)$ is said to be
 $G$-BiHom-associative if the following equation holds:
\begin{eqnarray*}
\sum_{\sigma\in G}(-1)^{\varepsilon(\sigma)}(-1)^{|\sigma(x_{1},x_{2},x_{3})|}
    \textbf{as}_{\alpha,\beta}\circ\sigma(\alpha^{-1}\beta^{2}(x_{1}),\beta(x_{2}),\alpha(x_{3}))=0,
\end{eqnarray*}
for all homogeneous elements $x_{1},x_{2},x_{3}\in V$.
\medskip

 \noindent{\bf Proposition 3.7.}
Let $G$ be a subgroup of the symmetric group $S_3$.
Then any $G$-BiHom-associative superalgebra $(V,\mu,\alpha,\beta)$ is BiHom-Lie admisssible.
\smallskip

{\bf Proof}
The BiHom-supersymmetry (2.6) follows straightaway from the definition.
Assume that $G$ is a subgroup of $S_3$. Then $S_3$ can be written as the disjoint union of the left cosets of $G$.
Say $S_3=\bigcup_{\sigma\in I}$, with $I\subseteq S_{3}$, and for any $\sigma,\sigma'\in I$,
$\sigma\neq\sigma'\in I \Rightarrow \sigma G\cap\sigma' G=\emptyset.$
It follows that
\begin{eqnarray*}
&&\sum_{\sigma\in S_{3}}(-1)^{\varepsilon(\sigma)}(-1)^{|\sigma(x_{1},x_{2},x_{3})|}
    \textbf{as}_{\alpha,\beta}\circ\sigma(\alpha^{-1}\beta^{2}(x_{1}),\beta(x_{2}),\alpha(x_{3}))\\
 &&~~~~=\sum_{\tau\in I}\sum_{\sigma\in \tau G}(-1)^{\varepsilon(\sigma)}(-1)^{|\sigma(x_{1},x_{2},x_{3})|}
    \textbf{as}_{\alpha,\beta}\circ\sigma(\alpha^{-1}\beta^{2}(x_{1}),\beta(x_{2}),\alpha(x_{3}))=0,
\end{eqnarray*}
for all homogeneous elements $x_{1},x_{2},x_{3}\in V$.
By Lemma 3.6,  $(V,\mu,\alpha,\beta)$ is a  BiHom-Lie admisssible superalgebra.
The proof is completed.
\hfill $\square$
\medskip

Now we provide a classification of the BiHom-Lie admisssible superalgebras via
$G$-BiHom-associative superalgebras. The subgroups of $S_{3}$ are
\begin{eqnarray*}
&&G_{1}=\{id\},~~G_{2}=\{id,\sigma_{1}\},~~G_{3}=\{id,\sigma_{2}\},\\
&&G_{4}=\{id,\sigma_{2}\sigma_{1}\sigma_{2}=(13)\},~~G_{5}=A_{3},~~G_{6}=S_{3},
\end{eqnarray*}
where $A_3$ is the alternating subgroup of $S_3$.
\smallskip

(1) The $G_{1}$-BiHom-associative superalgebras are the BiHom-associative superalgebras defined in Definition 2.1.

(2) The $G_{2}$-BiHom-associative superalgebras satisfy the condition:
\begin{eqnarray*}
&&\mu(\beta^{2}(x),\mu(\beta(y),\alpha(z)))-\mu(\mu(\alpha^{-1}\beta^{2}(x),\beta(y)),\alpha\beta(z))\\
&&~~=(-1)^{|x||y|}\{\mu(\alpha\beta(y),\mu(\alpha^{-1}\beta^{2}(x),\alpha(z)))-\mu(\mu(\beta(y),\alpha^{-1}\beta^{2}(x)),\alpha\beta(z))\}.
\end{eqnarray*}

(3) The $G_{3}$-BiHom-associative superalgebras satisfy the condition:
\begin{eqnarray*}
&&\mu(\beta^{2}(x),\mu(\beta(y),\alpha(z)))-\mu(\mu(\alpha^{-1}\beta^{2}(x),\beta(y)),\alpha\beta(z))\\
&&~~=(-1)^{|y||z|}\{\mu(\beta^{2}(x),\mu(\alpha(z),\alpha(y)))-\mu(\mu(\alpha^{-1}\beta^{2}(x),\alpha(z)),\beta^{2}(y))\}.
\end{eqnarray*}

(4) The $G_{4}$-BiHom-associative superalgebras satisfy the condition:
\begin{eqnarray*}
&&\mu(\beta^{2}(x),\mu(\beta(y),\alpha(z)))-\mu(\mu(\alpha^{-1}\beta^{2}(x),\beta(y)),\alpha\beta(z))\\
&&~~=(-1)^{|x||y|+|x||z|+|y||z|}\{\mu(\alpha^{2}(z),\mu(\beta(y),\alpha^{-1}\beta^{2}(x)))-\mu(\mu(\alpha(z),\beta(y)),\alpha^{-1}\beta^{3}(x))\}.
\end{eqnarray*}

(5) The $G_{5}$-BiHom-associative superalgebras satisfy the condition:
\begin{eqnarray*}
&&\mu(\beta^{2}(x),\mu(\beta(y),\alpha(z)))-\mu(\mu(\alpha^{-1}\beta^{2}(x),\beta(y)),\alpha\beta(z))\\
&&~~=-(-1)^{|x||y|+|x||z|}\{\mu(\alpha\beta(y),\mu(\alpha(z),\alpha^{-1}\beta^{2}(x)))-\mu(\mu(\beta(y),\alpha(z)),\alpha^{-1}\beta^{3}(x))\}\\
&&~~~~~-(-1)^{|y||z|+|x||z|}\{\mu(\alpha^{2}(z),\mu(\alpha^{-1}\beta^{2}(x),\beta(y)))-\mu(\mu(\alpha(z),\alpha^{-1}\beta^{2}(x)),\beta^{2}(y))\}.
\end{eqnarray*}

(5) The $G_{6}$-BiHom-associative superalgebras are the BiHom-Lie admisssible superalgebras.

\section{Derivations of BiHom-Lie superalgebras}
\def\theequation{\arabic{section}. \arabic{equation}}
\setcounter{equation} {0}

 In this section, we provide the notion of derivations of a BiHom-Lie superalgebra  $L$ and
 prove that the set of all derivations of $L$  has a natural BiHom-Lie superalgebra structure.
\medskip

Let $L=(L,[\cdot,\cdot],\alpha,\beta)$ be a BiHom-Lie  superalgebra. For any nonnegative integer $k$,
denote by $\alpha^{k}$ the $k$-times composition of $\alpha$, i.e.
\begin{eqnarray*}
\alpha^{k}=\alpha\circ\cdots\circ\alpha ~(k\mbox{-times}).
\end{eqnarray*}
In particular, $\alpha^{-1}=0, \alpha^{0}=id$ and $\alpha^{1}=\alpha$.
And similarly for the notion $\beta^{k}$.
\medskip

 \noindent{\bf Definition 4.1.}
For any integer $k\geq -1$, a homogeneous linear map $D: L\rightarrow L$ of degree $|D|$
is called a $\beta^{k}$-derivation of the BiHom-Lie  superalgebra  $(L,[\cdot,\cdot],\alpha,\beta)$
 if it satisfies
\begin{eqnarray}
&&D\circ\alpha=\alpha\circ D, D\circ\beta=\beta\circ D,\\
&&D[x,y]=[D(x), \beta^{k}(y)]+(-1)^{|x||D|}[\beta^{k}(x),D(y)],
\end{eqnarray}
for all homogeneous elements $x,y\in L$.
\medskip

We denote by $Der_{\beta^{k}}(L)=(Der_{\beta^{k}}(L))_{0}\oplus(Der_{\beta^{k}}(L))_{1}$
the set of $\beta^{k}$-derivation of the BiHom-Lie  superalgebra $(L,[\cdot,\cdot],\alpha,\beta)$, and
$Der(L)=\bigoplus_{k\geq-1}Der_{\beta^{k}}(L)$.
Define the endomorphisms $\widetilde{\alpha},\widetilde{\beta}$ on $Der(L)$ by
\begin{eqnarray*}
\widetilde{\alpha}(D)=\alpha\circ D,~\widetilde{\beta}(D)=\beta\circ D.
\end{eqnarray*}
For any $D,D'\in Der(L)$, define their commutator $[D,D']$ as follows:
\begin{eqnarray*}
[D,D']=D\circ D'-(-1)^{|D||D'|}D'\circ D.
\end{eqnarray*}

 \noindent{\bf Lemma 4.2.}
Let $(L,[\cdot,\cdot],\alpha,\beta)$ be a BiHom-Lie  superalgebra.
For any $D\in(Der_{\beta^{k}}(L))_{i},~D'\in(Der_{\beta^{s}}(L))_{j}$,
where $k+s\geq-1$ and $(i,j)\in Z_{2}^{2}$,
then $[D,D']\in(Der_{\beta^{k+s}}(L))_{j}$.
\smallskip

{\bf Proof}
For any $x,y\in L$, we have
\begin{eqnarray*}
&&[D,D']([x,y])\\
&&~~~~=(D\circ D'-(-1)^{|D||D'|}D'\circ D)([x,y])\\
&&~~~~=D([D'(x),\beta^{s}(y)]+(-1)^{|x||D|}[\beta^{s}(x),D'(y)])\\
&&~~~~~~~-(-1)^{|D||D'|}D'([D(x), \beta^{k}(y)]+(-1)^{|x||D|}[\beta^{k}(x),D(y)])\\
&&~~~~=[DD'(x), \beta^{s+k}(y)]+(-1)^{|D||D'(x)|}[D'(\beta^{k}(x)),D(\beta^{s}(y))]\\
&&~~~~~~~+(-1)^{|x||D'|}([D(\beta^{s}(x)),D(\beta^{k}(y))]
          +(-1)^{|x||D|}[\beta^{s+k}(x),DD'(y)])\\
&&~~~~~~~-(-1)^{|D||D'|}([D'D(x),\beta^{s+k}(y)]+(-1)^{|D'||D(x)|}[D(\beta^{s}(x)),D'(\beta^{k}(y))])\\
&&~~~~~~~-(-1)^{|D|(|D'|+|x|)}([D'(\beta^{k}(x)),D(\beta^{s}(y))]+(-1)^{|x||D'|}[\beta^{s+k}(x),D'D(y)])\\
&&~~~~=[DD'(x)-(-1)^{|D||D'|}D'D(x), \beta^{s+k}(y)]\\
&&~~~~~~~+(-1)^{|x|(|D|+|D'|)}[\beta^{s+k}(x),(DD'-(-1)^{|D||D'|}D'D)(y)]\\
&&~~~~=[[D,D'](x), \beta^{s+k}(y)]+(-1)^{|x|(|[D,D']|)}[\beta^{s+k}(x),[D,D'](y)].
\end{eqnarray*}
It is easy to check that $[D,D']\circ\alpha=\alpha\circ [D,D'], [D,D']\circ\beta=\beta\circ [D,D']$,
which leads to $[D,D']\in Der_{\alpha^{k+s}}(L)$.
$\hfill \Box$
\medskip

\noindent{\bf Remark 4.3.}
Obviously, we have
\begin{eqnarray*}
Der_{\beta^{-1}}(L)=\{D\in End (L)|D\circ\alpha=\alpha\circ D, D\circ\beta=\beta\circ D,D[x,y]=0, \forall x,y\in L\}.
\end{eqnarray*}
Thus for any $D,D'\in Der_{\beta^{-1}}(L)$, we have $[D,D']\in Der_{\beta^{-1}}(L).$
\medskip

 \noindent{\bf Proposition 4.4.}
Let $(L,[\cdot,\cdot],\alpha,\beta)$ be a BiHom-Lie  superalgebra.
Then $(Der(L),[\cdot,\cdot],\widetilde{\alpha},\widetilde{\beta})$ is  a BiHom-Lie  superalgebra.
\smallskip

{\bf Proof}
We prove that  the bracket $[\cdot,\cdot]$ on $Der(L)$ satisfies the conditions in Definition 2.3.
Let $D\in(Der_{\alpha^{k}}(L))_{i},~D'\in(Der_{\alpha^{s}}(L))_{j},~D''\in(Der_{\alpha^{t}}(L))_{l}$ and $x\in L$,
we have
\begin{eqnarray*}
(\widetilde{\alpha}\circ\widetilde{\beta})(D)=D\circ\alpha\circ\beta=D\circ\beta\circ\alpha=(\widetilde{\beta}\circ\widetilde{\alpha})(D).
\end{eqnarray*}
So Eq. (2.4) holds and similarly for Eq. (2.5). For Eq. (2.6), we have
\begin{eqnarray*}
&&[\widetilde{\beta}(D),\widetilde{\alpha}(D')]
=[D\circ\beta,D'\circ\alpha]\\
&&~~~~=(D\circ\beta)\circ(D'\circ\alpha)-(-1)^{|D||D'|}(D'\circ\alpha)\circ(D\circ\beta)\\
&&~~~~=(D\circ D'-(-1)^{|D||D'|}D'\circ D)\circ(\alpha\beta)\\
&&~~~~=-(-1)^{|D||D'|}(D'\circ D-(-1)^{|D||D'|}D\circ D')\circ(\alpha\beta)\\
&&~~~~=-(-1)^{|D||D'|}[\widetilde{\beta}(D'),\widetilde{\alpha}(D)].
\end{eqnarray*}
For Eq. (2.7), we calculate
\begin{eqnarray*}
&&(-1)^{|D||D''|}[\widetilde{\beta}^{2}(D),[\widetilde{\beta}(D'),\widetilde{\alpha}(D'')]
=(-1)^{|D||D''|}[D\circ\beta^{2},[D'\circ\beta, D''\circ\alpha]\\
&=&(-1)^{|D||D''|}[D\circ\beta^{2},(D'\circ D'')\circ(\beta\alpha)-(-1)^{|D'||D''|}(D''\circ D')\circ(\beta\alpha)]\\
&=&(-1)^{|D||D''|}\{(D\circ(D'\circ D''))-(-1)^{|D||D'|}((D'\circ D'')\circ D)\}\circ(\beta^{3}\alpha)\\
&&-(-1)^{|D''|(|D|+|D'|)}\{(D\circ(D''\circ D'))-(-1)^{|D'|(|D|+|D''|)}((D''\circ D')\circ D)\}\circ(\beta^{3}\alpha).
\end{eqnarray*}
Therefore, one can check that
$\circlearrowleft_{D,D',D''}(-1)^{|D||D''|}[\widetilde{\beta}^{2}(D),[\widetilde{\beta}(D'),\widetilde{\alpha}(D'')]=0$,
as desired. And this finishes the proof.
$\hfill \Box$
\medskip

For any homogeneous elements $a\in L$ satisfying $\alpha(a)=a=\beta(a)$,
define $ad_{k}(a)\in End(L)$ by
\begin{eqnarray*}
ad_{k}(a)(x)=[a,\beta^{k}(x)],~ \forall x\in L.
\end{eqnarray*}

\noindent{\bf Proposition 4.5.}
Let $(L,[\cdot,\cdot],\alpha,\beta)$ be a BiHom-Lie  superalgebra and $a$ an homogeneous element in $L$.
Assume that the structure maps $\alpha$ and $\beta$ are bijective,
then $ad_{k}(a)$ is an  $\beta^{k+1}$-derivation, which we call inner $\beta^{k+1}$-derivation.
\smallskip

{\bf Proof}
For any homogeneous elements $x,y\in L$, on the one hand we have
\begin{eqnarray*}
ad_{k}(a)[x,y]
&=&[a,\beta^{k}[x,y]]
=[\beta^{2}(a),[\beta^{k}(x), \beta^{k}(y)]]\\
&=&-(-1)^{|a||y|}(-1)^{|x||a|}[\beta^{k+1}(x),[\beta^{k+1}\alpha^{-1}(y),\alpha(a)]]\\
&&-(-1)^{|a||y|}(-1)^{|x||y|}[\beta^{k+2}\alpha^{-1}(y),[\beta(a),\alpha\beta^{k-1}(x)]]\\
&=&-(-1)^{|a||y|}(-1)^{|x||a|}[\beta^{k+1}(x),[\beta^{k+1}\alpha^{-1}(y),\alpha(a)]]\\
&&-(-1)^{|a||y|}(-1)^{|x||y|}[\beta^{k+2}\alpha^{-1}(y),[a,\alpha\beta^{k-1}(x)]].
\end{eqnarray*}
On the other hand, we have
\begin{eqnarray*}
[ad_{k}(a)(x), \beta^{k+1}(y)]
&=&[[a,\beta^{k}(x)], \beta^{k+1}(y)]
=[\beta[a,\beta^{k-1}(x)], \beta^{k+1}(y)]\\
&=&-(-1)^{(|a|+|x|)|y|}[\beta\beta^{k+1}\alpha^{-1}(y),\alpha[a,\beta^{k-1}(x)]]\\
&=&-(-1)^{(|a|+|x|)|y|}[\beta^{k+2}\alpha^{-1}(y),[a,\alpha\beta^{k-1}(x)]]
\end{eqnarray*}
and
\begin{eqnarray*}
[\beta^{k+1}(x),ad_{k}(a)(y)]
&=&[\beta^{k+1}(x),[a,\beta^{k}(y)]]\\
&=&[\beta^{k+1}(x),[\beta(a),\alpha\beta^{k}\alpha^{-1}(y)]]\\
&=&-(-1)^{|a||y|}[\beta^{k+1}(x),[\beta\beta^{k}(y),\alpha(a)]]\\
&=&-(-1)^{|a||y|}[\beta^{k+1}(x),[\beta^{k+1}\alpha^{-1}(y),\alpha(a)]].
\end{eqnarray*}
It follows that
\begin{eqnarray*}
ad_{k}(a)[x,y]
=[ad_{k}(a)(x), \beta^{k+1}(y)]+(-1)^{|x||a|}[\beta^{k+1}(x),ad_{k}(a)(y)],
\end{eqnarray*}
as desired. And this finishes the proof.
\hfill $\square$
\medskip

%

\begin{center}
 {\bf ACKNOWLEDGEMENT}
 \end{center}

  The paper is partially supported by the China Postdoctoral Science Foundation (Nos. 2015M571725 and 2015M580508),
 the Key University Science Research Project of Anhui Province (Nos. KJ2015A294, KJ2014A183 and KJ2016A545),
  the  NSF of Chuzhou University (Nos. 2014PY08 and 2015qd01) and the Program for Science and Technology Innovation Talents in Education Department of Guizhou Province(No. KY[2015]481).

\renewcommand{\refname}{REFERENCES}

\end{document}